\RequirePackage{fix-cm}
\documentclass[smallextended]{svjour3}       
\smartqed  
\usepackage{graphicx,enumerate}
\usepackage{amsmath,amssymb}
\newcommand{\R}{\mathbb{R}}
\newcommand{\N}{\mathbb{N}}
\newcommand{\C}{\mathbb{C}}
\numberwithin{theorem}{section}
\numberwithin{corollary}{section}
\numberwithin{lemma}{section}
\numberwithin{remark}{section}
\numberwithin{equation}{section}
\begin{document}

\title{$O_M(\R^n)$ as locally convex Orlicz space}

\author{Jan Kisyński}

\institute{
 \email{jan.kisynski@gmail.com}
}


\maketitle

\begin{abstract}
The original M.Valdivia proof of his theorem on representation of the space $O_M$ uses results of S.Rolewicz concerning metric linear spaces and results of A.Grothendieck from
his theory of topological tensor product. We present a more direct proof.
\keywords{Multiplication \and convolution \and the M.Valdivia representation of $O_M$}
\end{abstract}

\section{Introduction}
\label{sec:1}
Usually the L.Schwartz set $S(\R^n)$ of infinitely differentiable functions on $\R^n$ rapidly decreasing at infinity is treated as a topological space with distinguished topology. But in fact $S(\R^n)$ admits at least two important locally convex topologies. These topologies are determined by systems of seminorms $\{\rho_{\alpha,\beta};\alpha\in\N,\beta\in\N^n_0\}$ and $\{\tilde{\rho}_{\alpha,\beta};\alpha\in\N,\beta\in\N^n_0\}$ whose exact definition will be given in Section \ref{sec:2}.

Let $S(\R^n)$ be the L.Schwartz set of infinitely differentiable functions on $\R^n$ rapidly decreasing at infinity and let $S' (\R^n)$ be its dual set. In 1981 M.Valdivia has shown in \cite{12} that for every element $\Phi$ the space $O_M (\R^n)$ there is a distribution $T$ on $\R^n$ rapidly decreasing at infinity such that

$$\Phi\cdot\psi=T*\psi\quad \textnormal{for every}\quad \psi\in S(\R^n)$$
where $T*\psi:=T_{r_0}(\check{\psi})$ is a distinguished regularization of T.M.Valdivia used theory of linear metric spaces developed in the book \cite{6} of S.Rolewicz and used also the theory of topological tensor products of A.Grothendieck \cite{4}. Our aim is to prove the Valdivia theorem more directly, remaining in the framework of $(S(\R^n),\tilde{\rho})$ and $(S(\R^n),\tilde{\rho})'$. We use the operators $(m\cdot)$ of multiplication of elements of $(S(\R^n),\tilde{\rho})$ by multipliers $m$ of $(S(\R^n),\tilde{\rho}$), and we distinguish the operators $(m\cdot)$ from the multipliers themselves. Moreover, we use the elementary convolution $(S'(\R^n),\tilde{\rho})'\times(S(\R^n),\tilde{\rho}) \ni(T , \psi)\to T (\check{\psi})\in(S(\R^n),\tilde{\rho})$. Notice that elementary convolutions are deeply related to regularization, the last being described in Sect.4.10 of John
Horv\'{a}th book \cite{6}. Distributions $T\in(S(\R^n),\tilde{\rho})'$ rapidly decreasing at $\infty$ are in [3, Sect.125.0] characterized by their regularizations. The distinguished regularization $T_{r_0}$ occurring in the proof of \ref{thm41} has the value at $\tilde{\psi}\in(S(\R^n),\tilde{\rho})$ equal to $F(\psi)$, the Fourier image of $\psi$. This resembles the Fourier–Schwartz transformation.

\section{Locally convex topologies in the L.Schwartz set $S(\R^n)$}
\label{sec:2}
Except of the L.Schwartz set $S(\R^n)$ of interest is a larger set $\tilde{S}(\R^n)$ whose element are these functions $\psi\in C^\infty(\R^n;\C)$ for which $\lim\limits_{|x|\to\infty}[\partial^\beta\psi](x)=0$ for every $\beta\in\N_0^n$ in sense of the usual convergence, not necessarily rapid. The set $\tilde{S}(\R^n)$ admits the locally convex topology $\rho$
determined by system of seminorms $\{\rho_{\alpha,\beta};\alpha\in\N,\beta\in\N_0^n\}$ such that
\begin{align*}
[\rho_{\alpha,\beta}(\psi)](x)=\sup_{x\in\R^n}(1+|x|^2)^{-\frac{\alpha}{2}}|[\partial^\beta\psi](x)|\ \textnormal{for every}\ x\in\R^n\\
\textnormal{and every}\ \psi\in\tilde{S}(\R^n).
\end{align*}
The set $\tilde{S}(\R^n)$ admits also the locally convex topology $\tilde{\rho}$, determined by the system of seminorms $\{\tilde{\rho}_{\alpha,\beta};\alpha\in\N,\beta\in\N_0^n\}$ such that
$$\tilde{\rho}_{\alpha,\beta}(\psi)](x)=\int\limits_{\R^n}(1+|x|^2)^{-\frac{\alpha}{2}}|[\partial^\beta\psi](x)|dx\ \ \textnormal{for every}\ \ \psi\in\tilde{S}(\R^n).$$
where the integral over $\R^n$ is an absolutely convergent integral in sense of Lebesgue. Indeed, from [2, Lemma 4.15] it follows that there holds

\begin{theorem}
\label{thm21}
For every $\alpha\in\N, \beta\in\N_0^n$ and $\psi\in\tilde{S}(\R^n)$ there is $\mu\in[0,\infty]$ such that
$$\tilde{\rho}_{\alpha,\beta}(\psi)\leq\mu\rho_{\alpha,\beta}(\psi).$$
Theorem \ref{thm21} implies the
\end{theorem}

\begin{corollary}
\label{cor21}
In the set $\tilde{S}(\R^n)$ the topology $\tilde{\rho}$ is no stronger that the topology~$\rho$.
\end{corollary}

\section{Multipliers of $(S(\R^n),\tilde{\rho})$}
\label{sec:3}
The function $m\in C^\infty(\R;\C)$ is called multiplier of the space $(S(\R^n),\tilde{\rho})$, if $m\cdot(S(\R^n),\tilde{\rho})\subset(S(\R^n),\tilde{\rho})$. The set of all multipliers of $(S(\R^n),\tilde{\rho})$ is denoted by $m(\R^n)$ and constitutes a complex multiplication algebra. Let $O_M(\R^n)$ be the set of functions belonging to $C^\infty(\R^n;\C)$, growing at the most polynomially at $\infty$. The above assumption about growth means exactly that for every $\Phi\in O_M(\R^n)$ and every multiindex $\beta\in\N_0^n$ there is $\alpha\in\N$ such that $\lim\limits_{|x|\to\infty}\sup(1+|x|^2)^{-\frac{\alpha}{2}}|[\partial^\beta\Phi](x)|<\infty$, in usual sense, i.e. non necessarily in sense of rapid convergence.

\begin{theorem}[$m(\R^n)=O_M(\R^n)$.]
\label{thm31}
The inclusion $O_M(\R^n)\subset m(\R^n)$ is almost evident. Indeed, if $\Phi\in O_M(\R^n)$ and
$\psi \in S(\R^n)$, the the product $\Phi\cdot\psi$ is an infinitely differentiable complex valued function on $\R^n$ rapidly decreasing to zero at infinity, so that $\Phi\in m(\R^n)$. Proof of the inclusion $m(\R^n)\subset O_M(\R^n)$ is more complicated and is postponed to the Appendix.
\end{theorem}

\begin{theorem}
\label{thm32}
If $m\in m(\R^n)$, then the mapping $(S(\R^n),\tilde{\rho})\ni\psi\to m\cdot\psi\in(S(\R^n),\tilde{\rho})$ is continuous.
\end{theorem}

\begin{proof}
Since the locally convex space $(S(\R^n),\tilde{\rho})$ is metrizable, its topology is characterized whenever characterized are the convergent countable sequences $(\psi_k)_{k\in\N}$ of functions belonging to $(S(\R^n),\tilde{\rho})$. We can characterize the convergence of such sequences is as follows:
\end{proof}

\begin{lemma}
The countable sequence $(\psi_k)_{k\in\N}\subset(S(\R^n),\tilde{\rho})$ is convergent to zero if and only if for every $\beta\in\N_0^n$, the following two conditions are satisfied:
\begin{align}
\label{eq31}
&\lim_{k\to\infty}\partial^\beta\psi_k=0\ \textnormal{almost uniformly on}\ \R^n,\\
\label{eq32}
&\lim_{|x|\to\infty}|x|^\alpha|\partial^\beta\psi_k(x)|=0\ \textnormal{uniformly with respect to}\ \alpha\in\N.
\end{align}
The limits occurring in these conditions are understood in usual sense of Calculus, i.e. not in sense of rapid convergence.
\end{lemma}

\begin{proof}
The assertions (\ref{eq31}) and (\ref{eq32}) express in terms of the Calculus the fact that the sequence $(\psi_k)_{k\in\N}$ converges to zero in sense of $(S(\R^n),\tilde{\rho})$, i.e. that $\lim\limits_{k\to\infty}\tilde{\rho}_{\alpha,\beta}(\psi_k)=0$ for every $(\alpha,\beta)\in\N\times\N_0^n$ where $\tilde{\rho}_{\alpha,\beta}(\psi_k)=\int\limits_{\R^n}(1+|x|^2)^{\frac{\alpha}{2}}|(\partial^\beta\psi_k)(x)|dx$. Therefore the countable sequence $(\psi_k)_{k\in\N}$ converges to zero in the topology of $(S(\R^n),\tilde{\rho})$ if and only if both the conditions (\ref{eq31}) and (\ref{eq32}) are satisfied.
\end{proof}

\begin{proof}[of Theorem 3.2]
We have to prove that if $m\in m(\R^n)$ and both the conditions (\ref{eq31}) and (\ref{eq32}) are satisfied, then analogous condition, say (\ref{eq31})$^*$ and (\ref{eq32})$^*$, are satisfied for the sequence $(m\cdot \psi_k)_{k\in\N}$. Indeed, the Leibniz rule (see [6, Sect.2.5, Proposition 3]) says that
$$\partial^\alpha(m\cdot \psi_k)=\sum_{\beta\leq\alpha}\binom{\alpha}{\beta}(\partial^{\alpha-\beta}(m))\cdot(\partial\beta\psi_k).$$
Hence (\ref{eq31})$^*$ is a consequence of (\ref{eq31}) and (\ref{eq32})$^*$ is a consequence of (\ref{eq32}).
\end{proof}

\begin{remark}
\label{rem34}
If $m\in m(\R^n)$ is given, then the mapping $(S(\R^n),\tilde{\rho})\ni\psi\to m\cdot\psi\in (S (\R^n),\tilde{\rho})$ is continuous.
\end{remark}

\section{Appendix. Proof of the inclusion $m(\R^n)\subset O_M(\R^n)$ and proof of the M.Valdivia representation of the space $O_M (\R^n)$ by rapidly decreasing distributions}
We shall divide this Appendix onto several subsections.

\subsection{$O_M(\R^n)$ as a locally convex Orlicz space.}
Let the elements of $O_M(\R^n)$ be denoted by $\Phi$. Then
\begin{align*}
O_M(\R^n)=\{\Phi\in C^\infty(\R^n;\C):\ \textnormal{for every}\ n\in\N\ \textnormal{the integral}\\
\int\limits_{\R^n}(1+|x|^2)^{-\frac{n}{2}}|[\partial^\beta\Phi](x)|dx\ \textnormal{is convergent}\}.
\end{align*}
The above follows from the fact that elements of OM are operators of multiplication acting in the space $(S(\R^n),\tilde{\rho})$ which admits the locally convex topology determined by the system of seminorms $\{\tilde{\rho}_{\alpha,\beta};\alpha\in\N,\beta\in\N_0^n\}$ such that for every $\psi\in S(\R^n)$ the Lebesgue integral
$$\tilde{\rho}_{\alpha,\beta}(\psi)=\int\limits_{\R^n}(1+|x|^2)^{-\frac{\alpha}{2}}|[\partial^\beta\Phi](x)|dx\ \textnormal{is convergent}.$$
It follows that $O_M(\R^n)$ is uniquely determined by the system of spaces\linebreak $\{S_{-n}: n\in\N\}$ consisting of convergent integrals
\begin{equation}
\label{eq41}
\int\limits_{\R^n}(1+|x|^2)^{-\frac{n}{2}}\sup_{|\beta|\leq b_n}|[\partial^\beta\Phi](x)|dx
\end{equation}
where $b_n\in\N$ is fixed for every space $S_{-n}$. Each $S_{-n}$ is a locally convex space with topology determined by system of seminorms equal to integrals occurring in (\ref{eq41}). The system of spaces $\{S_{-n}: n\in\N\}$ can be treated as a single locally convex Orlicz space denoted by $\tilde{O}$ whose topology is determined by the system of seminorms $\{s_\beta: \beta\in\N_0^n\}$ such that $s_\beta(\Phi)=\int\limits_{\R^n}(1+|x|^2)^{-\frac{n}{2}}|[\partial^\beta\Phi](x)|dx$ for every $\Phi\in\tilde{O}$, and the set of convergence factors is equal to $\{(1+|x|^2)^{-\frac{n}{2}}: n\in\N\}$. In definition of this locally convex Orlicz space the above seminorms and the above convergence factors are equally important.

\subsection{The locally convex space $((m(\R^n)\cdot), op_m)$ of operators of multiplication}
For every $m \in m(\R^n)$ denote by $(m\cdot)$ the operator of multiplication the elements of $S(\R^n)$ by $m$. The set $(m(\R^n)\cdot)$ of multiplication operators admits the locally convex topology determined by the set of seminorms $\{s_{\alpha,\beta,\psi}; \alpha\in\N,\beta\in\N_0^n,\psi\in S(\R^n)\}$ such that
$$s_{\alpha,\beta,\psi}((m\cdot))=\tilde{\rho}_{\alpha,\beta}(m\cdot\psi)\ \textnormal{for every}\ (m(\R^n)\cdot).$$
The so defined locally convex topology in the set $(m(\R^n)\cdot)$ will be denoted by $op_m$. Existence of this topology is crucial for Subsection 4.5. Notice that for compact subsets of $\R^n$ the topology $op_m$ coincides with the topology introduced in [6, Sect.2.4, Example 11 and Sect.2.5, Example 8]. The elementary convolution $T*\psi$ where $\psi\in S(\R^n)$ and $T$ is a distribution on $(S (\R^n),\tilde{\rho})$ was analysed in [13, Sect.VI.3].

\subsection{The openworkness of the set $(m(\R^n)\cdot)$}
The set $(m(\R^n)\cdot)$ is openwork (in sense of this adjective explained in common words in [12, p.724]) because for every $k\in\N$ there exists in $(m(\R^n)\cdot)$ the inner automorphism ([2, Sect.190.D], [7, Sect.I.7])
$$A_k:(m(\R^n)\cdot)\ni(m\cdot)\to F_o^k\circ(m\cdot)\circ F_o^{-k}\in(m(\R^n)\cdot).$$
In the above formula $F_o$ denotes the Fourier transformation restricted to $S(\R^n)$. Let us stress that using the bracket-less symbol of superposition of operators we admit the rule that the operators act in order from right to left. For our proof of Theorem \ref{thm31} the openworkness of $(m(\R^n)\cdot)$ is only nonessential curiosity.

\subsection{Proof of the inclusion $m(\R^n)\subset O_M(\R^n)$}
\begin{theorem}
\label{thm41}
The inclusion
\begin{equation}
\label{eq42}
m(\R^n)\subset O_M(\R^n)
\end{equation}
is true
\end{theorem}

\begin{proof}
Usual multiplication of complex numbers implies that the inclusion (\ref{eq42}) is equivalent to the inclusion
\begin{equation}
\label{eq43}
m(\R^n)\cdot\psi\subset O_M(\R^n)\cdot\psi\ \textnormal{for every}\ \psi\in(S(\R^n),\tilde{\rho}).
\end{equation}
Left side of \ref{eq43} is equal to the set of values of the operators $(m(\R^n)\cdot)$ acting in $S(\R^n)$. Right side of \ref{eq43}  is equal to the set of values operators $\{\phi\in O_M(\R^n)\}$ which, by definition of $O_M(\R^n)$, is the maximal set of multiplication operators acting in $S(\R^n)$. \qed
\end{proof}

\subsection{Proof of M.Valdivia representation of the space $O_M(\R^n)$ by rapidly decreasing distributions}
We shall consider actions by multiplication. All the used properties of multiplication are
consequences of the usual properties of multiplication of complex numbers. Recall that a function $m\in C^\infty(\R^n;\C)$ is called a multiplier of the space $(S(\R^n),\tilde{\rho})$ if $m\cdot(S(\R^n),\tilde{\rho})\in(S(\R^n),\tilde{\rho})$. The above (very popular) definition concerns a single element of $C^\infty(\R^n;\C)$ and is not useful if
one wants to speak about action of a set of operators of multiplication in a set of functions.

We define the action of the set $(m(\R^n)\cdot)$ of operators of multiplication in the set
$m(\R^n)$ of multipliers of the space $(S(\R^n),\tilde{\rho})$ by the equality
\begin{equation}
\label{eq44}
(m(\R^n)\cdot)\cdot m(\R^n):= m(\R^n)\cdot(S(\R^n),\tilde{\rho}).
\end{equation}

Using theory of regularization of distributions exposed by John Horv\'{a}th in Section 4.10 of
his book [6], and using the definition of distribution on $\R^n$ rapidly decreasing at $\infty$ due in [3, Sect.125.0], one concludes that the representation theorem of M.Valdivia is equivalent to the following Theorem \ref{thm42}.

\begin{theorem}
\label{thm42}
The distribution $T\in S(\R^n)$ rapidly decreases at $\infty$ if and only if the set $ C^\infty(\R^n;\C)$ admits the locally convex topology determined by the system of seminorms $\{s_{\tau}:\tau\in\Xi\}$ such that if $\psi\in C^{\infty}(\R^n;C)$, then $s_{\tau}(\psi)\in F(\psi)$ for every $\tau\in\Xi$
\end{theorem}

\begin{proof}
The formulas occurring in Theorem \ref{thm42} mean that there is a unique regularization $\check{R}$ defined on $(S(\R^n),\tilde{\rho})'$ such that $\check{R}(T)=T(\psi)$ for every $T\in(S(\R^n),\tilde{\rho})'$ and every $\psi\in(S(\R^n),\tilde{\rho})$. Existence and uniqueness of such regularization are obvious.\qed
\end{proof}

\subsection{A remark concerning the Gagliardo-Nirenberg equality}
Consider the metrizable linear space $C_B\subset C^\infty(\R^n;\C)$ and its imbedding into Lebesgue space $L(\R^n)$. In their book \cite{2} R.A.Adams and J.J.F.Fournier observed that existence of the imbedding of $C_B$ into $L(\R^n)$ is equivalent to the Lemma 4.15 of \cite{2} which, by turns, is equivalent to the simple behaviour of $C_B$ described on p.79 of \cite{2}. This observation (as well as sole Lemma 4.15) is surprising. Also its connections with famous Gagliardo-Nirenberg equality are interesting. For that matter we shall limit ourselves to superficial presentation of the related literature.\medskip

The papers I and II of E.Gagliardo and III of L.Nirenberg are quoted in the books \cite{10} and \cite{11}.
\begin{enumerate}[I.]
\item E.Gagliardo, \emph{Caretterizioni della trace sulla frontira relative ad alcune di funzioni in n variabili},\\
Rend. Sem. Mat. Padova 27 (1957), 284-305.
\item E.Gagliardo, \emph{Proprietá di alcune classi di funzioni in più variabili},\\
Ricerche Mat. Napoli 7 (1958), 102-137.
\item L.Nirenberg, \emph{On elliptic partial differential equations},\\
Ann. della Scuole Normale Superiore, Pisa 13 (1959), 116-169.
\end{enumerate}
Russian Academy of Sciences published the translations of E.Gagliardo papers I and II in
Matematica 5:4 (1961), 87-116.\\[5pt]
Haim Brezis and Petru Mironescu published in HaL a survey of papers related to Gagliardo-
Nirenberg equality.


%
%



\end{document}